\documentclass[12pt]{article}

\usepackage{amssymb, amscd}
\usepackage{amsmath}
\usepackage{amsthm}
\usepackage{epsfig}
\usepackage{psfig}
\usepackage{color}

\newtheorem{dfn}{Definition}[section]

\newtheorem{thm}[dfn]{Theorem}
\newtheorem{lem}[dfn]{Lemma}
\newtheorem{rem}[dfn]{Remark}

\newtheorem{ass}[dfn]{Assumption}

\newtheorem{prop}[dfn]{Proposition}
\newtheorem{cor}[dfn]{Corollary}

\def\proof{\par\medskip\noindent{\it Proof: }}

\newcommand{\Isom}{\operatorname{Isom}}

\newcommand{\length}{\operatorname{length}}

\def\ulim{\mathop{\hbox{$\om$-lim}}}

\def\geo{\partial_{\infty}}

\def\A{{\mathbb A}}
\def\C{{\mathbb C}}
\def\F{{\mathbb F}}
\def\G{{\mathbb G}}
\def\H{\mathbb H}

\def\L{{\mathbb L}}
\def\N{{\mathbb N}}
\def\Q{{\mathbb Q}}
\def\R{{\mathbb R}}

\def\T{{\mathbb T}}

\def\Ga{\Gamma}
\def\eps{\epsilon}
\def\al{\alpha}
\def\be{\beta}
\def\ga{\gamma}
\def\om{\omega}
\def\La{\Lambda}

\def\del{\delta}
\def\si{\sigma}

\def\ol{\overline}
\def\ul{\underline}

\def\t{\tilde}
\def\acts{\curvearrowright}
\def\<{\langle}
\def\>{\rangle}

\begin{document}

\title{On sequences of finitely generated discrete groups}
\author{Michael Kapovich}
\date{\today}
\maketitle

\begin{abstract}
We consider sequences of  discrete subgroups $\Ga_i=\rho_i(\Ga)$ of a rank 1 Lie group $G$, 
with $\Ga$ finitely generated.  
We show that, for algebraically convergent sequences $(\Ga_i)$, 
unless $\Ga_i$'s are (eventually) elementary or contain normal finite subgroups 
of arbitrarily high order, their algebraic limit is a discrete nonelementary 
subgroup of $G$. In the case of divergent sequences $(\Ga_i)$ we show that the resulting 
action $\Ga\acts T$ on a real tree satisfies certain semistability condition, 
which generalizes the notion of stability introduced by Rips. We then verify that  
the group $\Ga$ splits as an amalgam or HNN extension of finitely generated groups, 
so that the edge group has an amenable image in $\Isom(T)$. 
\end{abstract}

\section{Introduction}

One of the basic questions about discrete subgroups of Lie groups is to understand limiting 
behavior of sequences of such groups. In this paper, we consider finitely generated 
discrete subgroups of a rank 1 semisimple Lie group $G$, i.e., generalizations 
of the classical Kleinian groups. Given a finitely generated 
group $\Ga$ and a sequence of  
subgroups $\Ga_i=\rho_i(\Ga)\subset G$, one says that this sequence {\em converges algebraically} 
to a subgroup $\Ga_\infty\subset G$ 
if the  sequence of homomorphisms $\rho_i: \Ga\to G$ 
converges (pointwise) to an epimorphism $\rho_\infty: \Ga\to \Ga_\infty\subset G$. 

More generally, one can consider algebraic convergence of $G$--equivalence classes 
of the representations $\rho_i$, 
where  $\rho_i$'s are replaced with their projections 
to the character variety ${\mathfrak X}(\Ga, G)$. Sequences $\rho_i$ which do not 
subconverge even in this sense, are called {\em divergent}. Every divergent sequence 
$(\rho_i)$ yields a nontrivial action  
$\Ga\acts T$ of the group $\Ga$ on a real tree $T$. One can regard such action 
as a generalization of the algebraic limit of the sequence $\Ga_i$. 

If the groups $\Ga_i$ are {\em discrete} and the representations $\rho_i$ are {\em faithful} 
then the limiting behavior is completely understood due to the following theorems:

\begin{thm}\label{1}
 Suppose that the group $\Ga$ is not virtually nilpotent 
and the sequence $\Ga_i$ converges algebraically to $\Ga_\infty$. 
Then  the algebraic limit $\Ga_\infty$ is discrete and $\rho_\infty$ is faithful.
\end{thm}

The above theorem is due to V.~Chuckrow \cite{Chuckrow}, 
N.~Wielenberg \cite{Wielenberg}, T.~Jorgenesen \cite{Jorgensen(1976)}, 
G.~Martin \cite{Martin(1993)}, in different degrees of generality, see  \cite{Martin(1993)} 
for the most general statement. 

\begin{rem}
Historically, it was H.~Poincar\'e \cite{Poincare} 
who first (unsuccessfully) tried to prove Theorem \ref{1} for Fuchsian subgroups of 
$SL(2,\R)$ as a part of his first attempt on proving the Uniformization 
Theorem (via the continuity method). 
\end{rem}

\begin{thm}
Suppose that the group $\Ga$ is not virtually nilpotent and 
the sequence $(\rho_i)$ is divergent (in the character variety). Then the limiting group 
action on the tree $\Ga\acts T$  is such that:

1. $\Ga\acts T$ is {\em small}, i.e. arc stabilizers are virtually solvable. 

2. The action $\Ga\acts T$ is {\em stable}. The group $\Ga$ splits as 
$\Ga=\Ga_1*_E \Ga_2$ or  $\Ga=\Ga_1*_E$ with the edge group $E$ amenable. 
\end{thm}

The first part is due to J.~Morgan and P.~Shalen \cite{Morgan-Shalen(1984)}, 
J.~Morgan \cite{Morgan(1986)}, M.~Bestvina \cite{Bestvina(1988)} 
and F.~Paulin \cite{Paulin(1988)} in the case when $G=SO(n,1)$. 
The proof in the case of other rank 1 Lie groups follows, for instance, by 
repeating the argument using the ultralimits which can be found in 
\cite[Chapter 10]{Kapovich2000}.  The second part, for finitely-presented groups, 
is mostly due to I.~Rips; see 
\cite{Rips-Sela(1994), Bestvina-Feighn(1995), GLP2, Paulin(1996), Kapovich2000} 
for the proofs. The theorem was recently extended to the case when $\Ga$ is merely 
finitely generated by V.~Guirardel \cite{Guirardel(2006)}.

\medskip 
The main goal of this paper is to analyze the case when the groups $\Ga_i$ are discrete 
but the representations $\rho_i$ are not necessarily faithful. As far as 
convergent sequences of discrete groups, the best one can hope for is to show that $\Ga_\infty$ 
is discrete and nonelementary, provided that the groups $\Ga_i$ are also discrete and 
nonelementary. This was proven by T.~Jorgensen and P.~Klein \cite{Jorgensen(1982)} 
in the case when $G=SL(2,\C)$ by methods specific to the 3-dimensional hyperbolic geometry. 
G.~Martin \cite{Martin(1989a)} observed that already for the hyperbolic 4-space, discreteness of 
$\Ga_\infty$ can fail. His example consisted of groups $\Ga_i=\Ga_i'\times \Phi_i$, where each 
$\Ga_i$ preserves a hyperbolic plane $\H^2\subset \H^4$ and the groups $\Phi_i$ are finite 
cyclic groups, so that the 
generators of $\Phi_i$ converge to a rotation of infinite order about $\H^2$. 
Martin proved in \cite{Martin(1989a)} for $G=SO(n,1)$ and in \cite{Martin(1993)} for 
isometry groups of negatively pinched Hadamard manifolds, 
that $\Ga_\infty$ is discrete and nonelementary provided that the groups $\Ga_i$ 
have {\em uniformly bounded torsion}. (See also \cite[Proposition 8.9]{Kapovich2000},
and \cite{Belegradek1998} for the proofs of discreteness of {\em geometric} limits, under the 
same assumption of uniformly bounded torsion and \cite{Wang-Yang} for another variation on 
the bounded torsion condition.) The uniform bound on torsion allows one to reduce  
the arguments to analyzing certain torsion-free elementary subgroups of $G$; 
such groups have the following property:

If $\La_1, \La_2$ are torsion-free discrete elementary subgroups of $G$, so that 
$\La_1\cap \La_2$ is nontrivial, then $\La_1, \La_2$ generate an elementary subgroup of $G$. 

It is easy to see that this property fails for subgroups with torsion and this is where  
the arguments of \cite{Martin(1989a), Martin(1993), Kapovich2000, Belegradek1998} 
break down in the presence of unbounded torsion. 

Our first result is

\begin{thm}\label{main1}
Suppose that $\Ga_\infty$ is an algebraic limit of a sequence of discrete nonelementary 
subgroups $\Ga_i\subset G$. Then: 

1. $\Ga_\infty$ is nonelementary. 

2. If $\Ga_\infty$ is nondiscrete, then for every sufficiently large $i$,  
each $\Ga_i$ preserves a proper symmetric subspace $X_i\subset X$. 
The kernel $\Phi_i$ of the restriction map $\Ga_i\to \Isom(X_i)$ is a finite
subgroup whose order $D_i$ diverges to infinity as $i\to\infty$. 

3. Every element $\ga$ of $\ker(\rho_\infty)$ either belongs to $\ker(\rho_i)$ for all 
sufficiently large $i$, or $\rho_i(\ga)\in \Phi_i$, where $\Phi_i$ is as in 2. 
\end{thm}

Therefore,  the example of G.~Martin described above is, in a sense, the only 
way the group $\Ga_\infty$ may fail to be discrete. 
(See remarks in the end of section 5 in \cite{Martin(1989a)}.) 
 In Corollary \ref{geo} we generalize Theorem \ref{main1} to geometric limits 
 of algebraically convergent sequences $\Ga_i$.

\medskip  
Our second result deals with the group actions on trees. Suppose that $\Ga$ is 
finitely generated, 
the groups $\Ga_i$ are discrete and the sequence $(\rho_i)$ diverges in the character variety. 
In general, there is no reason to expect  the action $\Ga\acts T$ to be {\em stable}. 
In Section \ref{semistability} we introduce the notion of {\em semistable} actions to remedy 
this problem. This notion requires stabilization not of sequences of arc stabilizers
$$
\Ga_{I_1}\subset \Ga_{I_2} \subset ...
$$
(as in the Rips' notion of stability) but stabilization of their {\em algebraic hulls}
\begin{equation}\label{a}
\A(\Ga_{I_1})\subset \A(\Ga_{I_2}) \subset ...
\end{equation}
which are certain solvable subgroups of $\Isom(T)$
canonically attached to $\Ga_{I_k}$. In the case at hand, the subgroups
$\A(\Ga_{I_k})$ are connected algebraic subgroups of a certain nonarchimedean
Lie group $\ul{G}(\F)$, for which $T$ is the Bruhat-Tits tree.
Stabilization of the sequence (\ref{a}) then comes from the fact that
the dimensions of the groups $\A(\Ga_{I_k})$ eventually stabilize.

\begin{rem}
M.~Dunwoody in his recent
preprint \cite{Dunwoody(2006)} proposed another way to eliminate the stability
 assumption for group actions on trees with {\em slender} arc stabilizers. 
 (A group is called {\em slender}
 if every subgroup is finitely generated.) However, both {\em slender} assumption is too
 restrictive (for instance, it forces the kernel of the action $\Ga\acts T$ to be 
 slender)  and the conclusion that Dunwoody obtains is not as strong as one would
 like.
\end{rem}

We then verify that semistability is sufficient for the Rips theory to work. 
As the result we obtain:

\begin{thm}
\label{main2}
Let $\rho: \Ga\acts T$ be the limiting action arising from a divergent sequence $(\rho_i)$. 
Then:

1. The action on $T$ of the image group $\bar\Ga:=\rho(\Ga)\subset \Isom(T)$ is small. 

2. The action $\rho: \Ga\acts T$ is  semistable.  

3. Assume that $\Ga$ is finitely-presented.  Then 
$\Ga$ splits as $\Ga=\Ga_1*_E \Ga_2$ or $\Ga=\Ga_1*_E$, where
$\rho(E)$ is a virtually solvable subgroup of $\Isom(T)$, and the groups 
$\Ga_1, \Ga_2, E$ are finitely generated.  
\end{thm}

In Propositions \ref{ker} and \ref{kern}, we also describe the kernel of the action $\Ga\acts T$. 

\medskip 
The key technical ingredient in the proof of  Theorems \ref{main1} 
and \ref{main2} is the definition of the {\em algebraic hull} $\A(\La)$ for amenable 
subgroups $\La\subset \G=\ul{G}(\L)$, where $\ul{G}$ is a reductive algebraic group 
and $\L$ is a field of cardinality continuum and zero characteristic. 
The group $\A(\La)$ is a (Zariski) connected algebraic solvable subgroup of 
${\G}$ so that the intersection 
$$
A(\La):=\A(\La)\cap \La
$$
is a subgroup of uniformly bounded index in $\La$. (The bound depends only on ${\G}$.)

\medskip
The results of this paper probably generalize to sequences of isometric 
group actions $\rho_i: \Ga\acts X_i$, where $X_i$ are Hadamard manifolds of 
fixed dimension with fixed pinching constants. However, at the moment, 
I am not sure how to establish such a generalization, 
as the concept of {\em algebraic hull} is missing in this context.

\medskip
{\bf Acknowledgements.} During this work the author was
partially supported by the NSF grants DMS-04-05180 and DMS-05-54349.
Most of this work was done when the author
was visiting the Max Plank Institute for Mathematics in Sciences located
in Leipzig. This work was motivated by conversations with Lena Klimenko and with 
Steve Boyer and Michel Boileau  in Summer of 2007: They
independently asked me to prove certain versions of 
Theorems \ref{main1} and \ref{main2} respectively. I am grateful to Igor Belegradek for
several useful comments. 

\section{Preliminaries}

For a $CAT(0)$ space $X$ we let $\geo X$ denote its visual boundary. We let $\ol{xy}\subset X$ 
denote the geodesic segment between points $x, y\in X$.  

Let $G$ be a  Lie group. Then there exists a neighborhood $U$ of $1$ in $G$, called 
{\em Zassenhaus neighborhood}, so that for every discrete subgroup $\Ga\subset G$, generated by 
elements from $U$, it follows that $\Ga$ is nilpotent. See e.g. \cite{Raghunathan}. 

Let $X$ be a {\em negatively pinched} Hadamard manifold, i.e. a complete simply-connected 
Riemannian manifold whose sectional curvature is bounded by two negative constants: 
$$
-a^2\le K_X \le -1.
$$
Then there exists a constant $\mu$ (called {\em Margulis constant}) 
which depends only on $a$ and the dimension of $X$ so that the following holds. For every point $x\in X$
and a discrete subgroup $\Ga$ of $\Isom(X)$ generated by elements $\ga_i$ so that 
$d(x, \ga_i(x))\le \mu$, it follows that the group $\Ga$ is virtually nilpotent. 
See e.g. \cite{Ballmann-Gromov-Schroeder(1985)}. 

Let $X$ be a negatively pinched Hadamard manifold. A discrete subgroup $\La\subset \Isom(X)$ 
is called {\em elementary} if one of the following equivalent conditions is satisfied:

a. $\La$ is amenable.

b. $\La$ contains no free nonabelian subgroups. 

c. $\La$ is virtually nilpotent.

d. $\La$ either fixes a point in $\bar{X}=X \cup \geo  X$, or preserves a geodesic
in $X$.

We refer the reader to \cite{Belegradek-Kapovitch} for a detailed description of the structure 
and quotient spaces $X/\La$ for such subgroups.

\section{Amenable subgroups of algebraic groups}

Let $G=\ul{G}(\C)$ be a connected reductive complex-algebraic Lie group, where $\ul{G}$
is defined over $\Q$. We will consider amenable subgroups
$\La\subset G$.

\begin{thm}\label{hull}
There exists a number $d=d(G)$ so that the following holds. For every
amenable subgroup $\La\subset G$ there exists a canonical (Zariski)
connected solvable algebraic subgroup $\A(\La)\subset G$ (the {\em algebraic hull}
of $\La$) so that:

1. $|A(\La): \La|\le d(G)$, where $A(\La):=\A(\La)\cap \La$.

2. $\A(\La)$ is canonical in the following sense:

a. If $\La_1\subset \La_2$ then $\A(\La_1)\subset \A(\La_2)$.

b. For every automorphism $\phi$ of $G$ (either algebraic, or coming from
$Gal(\C)$),
$$
\phi \A(\La) = \A(\phi \La ).$$

c. For every $g\in G$, if $g \A(\La) g^{-1}\subset \A(\La)$, then
$g \A(\La) g^{-1}= \A(\La)$.

\end{thm}
\proof We first, for the sake of being concrete, define $\A(\La)$ in the case
$G=SL(2, \C)$. Amenable  subgroups $\La\subset G$  are classified as follows:

1. $\La$ is finite of order $\le 120$; we then let $\A(\La):=\{1\}$.

2. $\La$ is finite of order $>120$; it then contains an abelian subgroup $A(\La)$
of index $\le 2$. The group $A(\La)$ is contained in a unique maximal torus
$\C^\times\cong \T\subset G$.  We then let $\A(\La):=\T$.

We now assume that $\La$ is infinite.

3. The subgroup $\La$ is diagonalizable. Then it
is contained in a unique maximal torus $\T\subset G$ (which has to be unique). Set 
$\A(\La):=\T$.

4. $\La$ is contained in the index 2 extension of a maximal torus $\T\subset G$.
We then let $\A(\La):=\T$.

5. $\La$ has a unique fixed point $\xi$ in $S^2=\geo  \H^3$. We then let $\A(\La)$
be the full stabilizer of $\xi$ in $G$. Up to conjugation, this
group consists of upper--triangular matrices in $G$ and is, therefore, solvable.

\medskip
We now discuss the general case.

1. Let $\Phi\subset G$ be a finite subgroup. Then (up to conjugation)
$\Phi$ is contained in the maximal compact subgroup $K\subset G$.
According to the Jordan Theorem, see e.g. \cite[Theorem 8.29]{Raghunathan},
there exists a canonical torus $\T=\T(\Phi)\subset K$, so that the
abelian subgroup $T(\Phi)=\Phi\cap \T(\Phi)$ has index $\le a(G)$ in $\Phi$.
We then let $\A(\Phi)\subset G$ be the complexification of the torus $\T$.

\medskip
Let $\La\subset G$ be an infinite amenable subgroup. Then, by the Tits alternative,
the Zariski closure $H:=\bar\La\subset G$ has to be virtually solvable.

2. Suppose that $H$ is an infinite reductive subgroup of $G$, i.e., its
Zariski component of the identity is a nontrivial
torus $\T_H\subset H$.  (This torus is not necessarily maximal.) Since $H$ has only
finitely many components, the quotient $\La/(\La\cap \T_H)$ is finite. The 
torus $\T_H$ is contained in the unique smallest torus $\T$ which is the 
intersection of maximal tori in $G$. The torus $\T$ corresponds under the 
exponential map to a face of a Weyl chamber of $G$. Therefore, the number 
of  conjugacy classes of such tori $\T \subset G$ is finite.

The group $\La$ is contained in $N(\T)$, the normalizer of the torus $\T$ in $G$.
Let $Z(\T)$ denote the centralizer of the torus $\T$ in $G$. Recall that
$$
N(\T)/Z(\T)=W_\T
$$
is the Weyl group associated with the torus $\T$. Hence, its order is bounded
from above by a constant $b=b(G)$. Therefore, $\La$ contains a subgroup $\La'$
of index $\le b$, so that $\La'\subset Z(\T)$. The quotient
$\La'/(\La'\cap \T)$ is a finite subgroup $\Phi$ of the Lie group
$Q:=Z(\T)/\T$. Since the number of conjugacy classes of the tori $\T \subset G$ is
finite, the number of components of $Q$ is bounded from above by some
$c=c(G)$. Therefore (by Case 1), there exists a canonical torus
$\A(\Phi)\subset Q$ so that
$$
|\Phi: A(\Phi)|\le a(Q),
$$
where $A(\Phi)=\A(\Phi)\cap \Phi$. The sequence
$$
1\to \T\to Z(\T)\to Q \to 1
$$
splits and we obtain
$$
\La'':=\La'\cap \T \times \A(\Phi), |\La: \La''|\le d:=a(Q)b(G).
$$
We then set $\A(\La):= \T \times \A(\Phi)$.

\medskip
3. Suppose that $H=\bar\La$ is a non-reductive subgroup of $G$.
Let $U\subset H$ be the unipotent radical of $H$, so $M:=H/U$ is reductive.
Let $\pi: H\to M$ denote the canonical projection.
The subgroup $U$ is solvable and is canonically defined. The Levi subgroup $M\subset G$ is
again algebraic. Therefore, we apply Case 2 to the subgroup $\pi(\La)\subset M$.
Then we set
$$
\A(\La):= \pi^{-1}(\A(\pi(\La)).
$$
Since $\A(\La)=U\cdot \A(\pi(\La))$ is the semidirect product of two solvable groups,
it is solvable itself.

Lastly, we verify the fact that $\A(\La)$ is canonical. Property (a) and invariance 
under algebraic automorphisms of $G$ follow from 
the construction. Consider invariance under the automorphisms $\phi$ of $G$ induced by 
$\si^{-1}\in Gal(\C)$. It suffices to treat the case when 
$G$ is an affine algebraic group (i.e., $GL(n, \C)$). 
Let $E\subset G$ be a subset and $f\in \C[G]$ 
be a polynomial function vanishing on $E$. Then $f^{\si}$ vanishes on $\phi(E)$. 
Moreover, if the ideal generated by the functions $f$ determines an algebraic 
subgroup $H$ of $G$, the same is true for the ideal generated by the functions 
$f^{\si}$. The subgroup $H$ is solvable and connected iff the corresponding 
subgroup $H^{\si}$ is. 

To check property (c) note that $g \A(\La) g^{-1}\subset \A(\La)$ implies that 
the above groups have the same Lie algebra. Then the equality follows from the 
connectedness of  $\A(\La)$. \qed

\begin{cor}\label{c:hull}
Let $\F$ be a field of characteristic zero and cardinality continuum,
$G=\ul{G}(\F)$ be an algebraic group. Then there exists a constant $d=d(G)$ so
that the following holds. Let $\La\subset G$ be an amenable subgroup.
Then there exists a canonical (Zariski) connected solvable algebraic subgroup
$\A(\La)\subset G$ so that: $|A(\La): \La|\le d(G)$, where $A(\La):=\A(\La)\cap \La$.
\end{cor}
\proof Let $\bar\F$ denote the algebraic closure of $\F$. Then
$\F\cong \C$  since both are extensions of $\Q$, algebraically closed and have the
same cardinality. Therefore, we may regard $\La$ as a subgroup of $\ul{G}(\C)$.
Let $\bar\A(\La)$ denote the algebraic hull of $\La\subset \ul{G}(\C)$.
Then, since $\bar\A(\La)$ is canonical, for every $\si\in Gal(\C/\F)$, we have
$$
\si(\bar\A(\La))=\bar\A(\La)
$$
We set
$$
\A(\La):=\bar\A(\La)\cap G.
$$
Then $\A(\La)$ is again solvable and Zariski connected. The rest of the properties
follow from Theorem \ref{hull}. \qed

\medskip
We will apply the above corollary in the following cases:
$\F=\R$ and $G$ is a real Lie group of rank 1; $\F$ is a complete
nonarchimedean valued field of zero characteristic and $G$ has rank 1.

\section{Algebraic limits of sequences of discrete groups}\label{limits}

In this section we prove Theorem \ref{main1}. Let $X$ be a negatively curved symmetric space; 
its isometry group is isomorphic to a (real) rank 1 
algebraic group $G$ defined over $\Q$.  For instance, the reader can think of $G=SO(n,1)$ and 
$X=\H^n$. Let $\rho_i: \Ga\to G$ be a sequence of discrete (but not necessarily faithful)
representations of a finitely generated group $\Ga$. We let $\Ga_i$ denote
the image of $\rho_i$. Suppose that $\lim_i \rho_i=\rho_\infty$ and 
$\Ga_\infty=\rho_\infty(\Ga)$ is the algebraic limit of the sequence $(\Ga_i)$. 
In the ``generic'' case, the group $\Ga_{\infty}:=\rho_\infty(\Ga)$ is a 
discrete nonelementary subgroup of $G$. The theorem below describes what 
happens in the exceptional cases.

\begin{thm}\label{c}
1. If $\Ga_{\infty}$ is discrete and elementary, then
for sufficiently large $i$, each $\Ga_i$ is elementary.

2. If $\Ga_\infty$ is nondiscrete, then either:

a. For every sufficiently large $i$, each $\Ga_i$ is elementary, or

b. For every sufficiently large $i$,  each $\Ga_i$ preserves a
proper symmetric subspace $X_i\subset X$. The kernel $\Phi_i$ of the restriction map
$\Ga_i\to \Isom(X_i)$ is a finite
subgroup whose order $D_i$ diverges to infinity as $i\to\infty$.
\end{thm}
\proof Let $U\subset G$ denote the Zassenhaus neighborhood of $1\in G$. 
Let $g_1,...,g_m$ denote the generators of $\Ga$.  We can assume that
$\Ga$ is free on the generators $g_1,...,g_m$. We will need

\begin{lem}\label{kernel}
Let $\ga\in \ker(\rho_\infty)$. Then for all but finitely many $i$ either 

(a) $\rho_i(\ga)=1$, or

(b) $\Ga_i$ is elementary, or

(c) $\Ga_i$ preserves a proper symmetric subspace $X_i\subset X$, which is fixed pointwise 
by $\rho_i(\ga)$. 
\end{lem}
\proof We assume that (a) does not occur. 
Let $K\subset \Ga$ denote the normal closure of $\{\ga\}$.
Exhaust $K$ by finitely generated subgroups
$$
K_1\subset K_2\subset ...
$$
so that
\begin{equation}\label{semicon}
g_j K_l g_j^{-1}\subset K_{l+1}, \quad \forall l, \forall j=1,...,m.
\end{equation}
Without loss of generality, we may assume that $\ga\in K_1$. 
It is standard that if $(h_i)$ is a sequence of nontrivial elements 
in a Lie group converging to $1$, then the orders of $h_i$ 
(regarded as elements of $\N\cup \{\infty\}$) converge to infinity.   

Therefore, since $\rho_i(\ga)\ne 1$ but $\lim_i  \rho_i(\ga)=1$, 
the order of $\rho_i(\ga)$ diverges to infinity
as $i\to \infty$ for each $j=1,...,s$. It follows that
the order of $\rho_i(K_1)$ diverges to infinity as $i\to\infty$.
In particular, without loss of generality, we may assume that for each $i$,
the hull $\A(\rho_i(K_1))$ is a nontrivial connected solvable subgroup of $G$.

For every $g\in K$, there exists $i_g$ so that for all $i\ge i_g$,
$\rho_i(g)\in U$. Therefore, without loss of generality, we may assume that
for all $i$, the groups
$$
\rho_i(K_l), l=1,..., D=\dim(G),
$$
are elementary, where $\dim(G)$ is the dimension of $G$. Hence,
for each $i$, there exists a pair of groups $\A(\rho_i(K_l)), \A(\rho_i(K_{l+1}))$ 
(for some $0\le l\le D-1$ depending on $i$) 
which have the same dimension, and, hence, are equal. These groups are necessarily
nontrivial.

Since $\A_{il}=\A(\rho_i(K_l))$ is canonical, in view of (\ref{semicon}) 
we obtain
\begin{equation}\label{eq}
\rho_i(g_j) \A_{il} \rho_i(g_j)^{-1}=\A_{i(l+1)}= \A_{il}, \forall j=1,...,m. 
\end{equation}

If the group $\A_{il}$ is noncompact, then it either has a unique fixed point
in $\geo X$ or an invariant geodesic. This point or a geodesic are invariant 
under $\Ga_i$ according to (\ref{eq}). Therefore, it follows that $\Ga_i$ is elementary 
in this case.

We next assume that $\A_{il}$ is compact for each $i, l$. 
By (\ref{eq}), the group $\Ga_i$ preserves the 
fixed--point set $X_i\subset X$ of $\A_{il}$, which is a 
symmetric subspace in $X$. Since
$$
|\rho_i(K_1): \A_{i1}\cap \rho_i(K_1)|\le d(G),
$$
we have
$$
\rho_i(\ga)^q \in \A_{i1} 
$$
for some $1\le q\le  d(G)$. Hence, $\rho_i(\ga)|X_i$ is an element of order $\le
d(G)$ of $\Isom(X_i)$.  Since $\rho_i(\ga)|X_i$ converge to $1$ as 
$i\to\infty$, it follows that $\rho_i(\ga)$  restrict trivially
to $X_i$ for all sufficiently large $i$, $j=1,...,k$. 
Therefore, $X_i$ is a proper subspace in $X$ invariant under $\Ga_i$. \qed 

\medskip 
We now continue with the proof of Theorem \ref{c}. 

1. Suppose that $\Ga_{\infty}$ is discrete and elementary. Then $\Ga_{\infty}$
is a lattice in a nilpotent Lie group with finitely many components.
In particular, $\Ga_{\infty}$ is finitely-presented.
It therefore has the presentation
$$
\< g_1,..., g_m| R_1,..., R_k\>
$$
where $R_1,...,R_k$ are words in the generators $g_1,...,g_m$. Since $\Ga$ is free, we
can regard these words as elements of $\Ga$. By Lemma \ref{kernel}, 
for all sufficiently large $i$ one of the following holds:

a. The group $\Ga_i$ is elementary.

b. $X$ contains a symmetric subspace $X_i$ invariant under $\Ga$, so that each 
$\rho_i(R_j), j=1,...,k$ restricts trivially to $X_i$. Therefore, 
$R_1,...,R_k$ belong to the kernel of the restriction homomorphism
$$
\Ga\to \Ga_i\to \Isom(X_i). 
$$
Therefore the homomorphism $\Ga\to \Isom(X_i)$ factors through $\Ga\to \Ga_\infty$. 
Thus, its image is an amenable group. Since the kernel of $\Ga_i\to \Isom(X_i)$ is amenable, 
it follows that $\Ga_i$ is itself amenable and, hence, elementary.

\medskip
Case 2. Suppose that $\Ga_{\infty}$ is nondiscrete. Our arguments are somewhat similar 
to the Case 1. 

Let $\bar\Ga_\infty$ denote the closure of $\Ga$ in $G$ with respect to the 
classical topology. Then the identity component $\bar\Ga_\infty^0$ of this group 
is a nontrivial nilpotent group, see e.g. \cite[Proposition 8.9]{Kapovich2000} 
or \cite[Lemma 8.8]{Belegradek1998}. In any case, $\Ga_\infty$ contains 
nontrivial elements $\ga=\rho_\infty(g)$ arbitrarily close to $1$.
As before, the order of such $\ga$ necessarily goes to infinity as $\ga$ 
approaches $1$.

Let $V$ be a neighborhood of $1$ in $G$ whose closure is contained
in the Zassenhaus neighborhood $U$.
By choosing $\ga$ sufficiently close to $1$, we obtain:
$$
\ga, \rho_\infty(g_j) \ga \rho_\infty(g_j)^{-1},..., \rho_\infty(g_j)^D \ga
\rho_\infty(g_j)^{-D} \in V, \quad j=1,...,m, 
$$
where $D$ can be taken as large as we like. Consider
the subgroups
$$
K_s:=\<g_j^t g g_{j}^{-t}, j=1,...,m, t=0,...,s\>\subset \Ga
$$
for $s=0,...,D$. Then,
$$
K_0\subset K_1\subset ... \subset K_D
$$
and
$$
g_j^t K_s g_{j}^{-t} \subset K_{s+1}, \forall j=1,...,m, s=0,...,D-1.
$$

As before, we choose $D$ equal the dimension of $G$. By considering sufficiently
large $i$ we can assume that
$$
\rho_i(g_j^s g g_{j}^{-s})\in U, j=1,...,m, s=0,...,D.
$$

Therefore, the subgroups $\La_{is}=\rho_i(K_s)$ generated by the above elements of $\Ga_i$, 
are elementary for $s=0,...,D$. Since $\ga$ can be taken to have arbitrarily high 
(possibly infinite) order, we can assume that the algebraic hull
$\A(\La_{is})$ is nontrivial for each $i$ and $s$.

We now repeat the arguments from the proof in Case 1. For each $i$,
there exists $0\le s<D$ so that
$$
\A(\La_{is})=\A(\La_{i(s+1)}).
$$
Therefore,
\begin{equation}\label{eq2}
\rho_i(g_j) \A(\La_{is}) \rho_i(g_j)^{-1} = \A(\La_{is}), j=1,...,m.
\end{equation}
If $\A(\La_{is})$ is noncompact, it follows from (\ref{eq2}) 
that $\Ga_i$ is elementary, which contradicts our assumptions. 
Therefore $\A(\La_{is})$ is compact (a torus in $G$); this subgroup 
fixes (pointwise) a proper symmetric subspace $X_i\subset X$. 
According to (\ref{eq2}), this
subspace is invariant under the group $\Ga_i$. The kernel $\Phi_i$ of the restriction
homomorphism $\Ga_i\to \Isom(X_i)$ contains the abelian
subgroup $A(\La_{is})=\A(\La_{is})\cap \La_{is}$. By construction, the order of
$A(\La_{is})$ diverges to infinity as $i\to\infty$. Therefore, the order $D_i$ of
$\Phi_i$ also  diverges to infinity as $i\to\infty$.
\qed

\begin{cor}
Suppose that $G=PSL(2,\C)$ and, hence, $X=\H^3$. Then:

1. Either $\Ga_{\infty}:=\rho_\infty(\Ga)$ is discrete and nonelementary,
or

2. For each sufficiently large $i$, the group $\Ga_i$ is elementary.
\end{cor}
\proof It suffices to analyze Case 2b of the above theorem.
Then $\Ga_i$ contains a nontrivial finite normal subgroup $\Phi_i$ of
rotations about a symmetric subspace $X_i\subset \H^3$;
this subspace is either a point or a geodesic. In either  case,
$\Ga_i$ is elementary. \qed

\begin{cor}\label{geo}
Suppose that $\Ga$ is a finitely generated group, 
homomorphisms $\rho_i: \Ga\to \Ga_i=\rho_i(\Ga)\subset 
G$ converge to $\rho_\infty: \Ga\to \Ga_\infty =\rho_\infty(\Ga)\subset G$ 
and the groups $\Ga_i$ are  discrete and nonelementary. 
Let $\Ga_\infty^{geo}\subset G$ be 
the geometric limit of the sequence of groups $\Ga_i$. Then: 

1. $\Ga_\infty^{geo}$ is nonelementary. 

2. If $\Ga_\infty^{geo}$ is nondiscrete, then   
each $\Ga_i$ contains a finite normal subgroup $\Phi_i$, whose order diverges 
to infinity as $i\to\infty$.
\end{cor}
\proof Recall that $\Ga_\infty\subset \Ga_\infty^{geo}$ (see e.g. \cite{Kapovich2000}). 
Since $\Ga_\infty$ is nonelementary 
by Theorem \ref{c}, it follows that $\Ga_\infty^{geo}$ is 
nonelementary as well. To prove Part 2, 
we modify Part 2 of the proof of Theorem \ref{c} as follows. 
Consider an  element $\ga\in \Ga_\infty^{geo}\setminus \{1\}$ 
sufficiently close to $1\in G$. Instead of using a fixed element $g\in \Ga$ 
so that $\rho_\infty(g)=\ga$, we consider a sequence $h_i\in \Ga$ so that
$$
\lim_{i\to\infty} \rho_i(h_i)=\ga. 
$$
Instead of the subgroups $K_s\subset \Ga$ we use  
$$
K_{s,i}:=\<g_j^t h_i g_{j}^{-t}, j=1,...,m, t=0,...,s\>\subset \Ga. 
$$
With these modifications, the proof of 
Part 2 of  Theorem \ref{c} goes through in the context of the geometric limit.  
\qed 

\section{Small actions}\label{small}

In this section we prove the first assertion of Theorem \ref{main2}.

Let $\rho: \Ga\acts T$ be an isometric action of a group $\Ga$ on a metric tree $T$. 
Let $\bar\Ga:=\rho(\Ga)\subset \Isom(T)$ denote the image of $\Ga$ in $\Isom(T)$. 
Given an axial isometry $g\in \Ga$, let $A_g$ denote the axis of $g$ and
$\ell(g)$ the translation length of $g$. Recall that the action $\bar\Ga\acts T$
is called {\em nontrivial} if $\Ga$ does not have a global fixed point.
This action is called {\em small} if the arc stabilizers are 
amenable.

Suppose that $(X,d)$ is a  negatively pinched simply-connected complete
Riemannian manifold and $\Ga$ is a finitely--generated group with the
generating set $\{g_1,...,g_m\}$.
Given a representation $\rho: \Ga\to \Isom(X)$, define
$$
b_x(\rho):= \max (d(\rho(g_i)(x), x), i=1,...,m),
$$
$$
b(\rho):= \inf_{x\in X} b_x(\rho).
$$

Then a sequence of representations $\rho_i: \Ga\to \Isom(X)$ is  
{\em divergent} if and only if 
$$
\lim_i b(\rho_i)=\infty.
$$
Indeed, if there is a subsequence $(\rho_{i_j})$ so that $b(\rho_{i_j})\le C$, 
then we can conjugate $\rho_{i_j}$ by the elements $h_{i_j}\in G$ which move $x_{i_j}$ 
to a base-point $o\in X$. Since $G$ is locally, compact, it follows that the new 
sequence 
$$
Ad(h_{i_j}) \rho_{i_j}
$$
converges in $Hom(\Ga, G)$. 

Let $\om$ be a nonprincipal ultrafilter on $\N$.
We recall that a divergent sequence yields a  nontrivial isometric action
$\rho_\om: \Ga\acts T$ of $\Ga$ on a metric tree $T$, well-defined up to
scaling (given the choice of $\om$). The tree $T$ is the $\om$-ultralimit
of the sequence of pointed metric spaces
$$
(X, \frac{d}{b(\rho_i)}, o_i)
$$
where $o_i\in X$ is the point {\em nearly}  realizing $b(\rho_i)$, i.e.,
$$
|b(\rho_i)-b_{o_i}(\rho_i)|\le 1.
$$
See e.g. \cite{Kapovich2000,Kapovich-Leeb(1995)} for the details.

We now assume that $X$ is a symmetric space, i.e. its isometry group is a rank 1
algebraic group $G$. 

\medskip
The following theorem is standard in the case of sequences of discrete
{\em and faithful} representations $(\rho_i)$:

\begin{thm}\label{smallness}
Let $\rho_i: \Ga\to \Isom(X)$ be a divergent sequence of representations with
discrete images. Let $\rho_\om: \Ga\acts T$ denote the limiting action on a tree 
and $\bar\Ga:=\rho_\om(\Ga)\subset \Isom(T)$. 
Then the action $\bar\Ga\acts T$ is  small.
\end{thm}
\proof Our proof repeats the arguments of the proof of Theorem 10.24 in
\cite{Kapovich2000} with minor modifications. 
Let $\mu>0$ denote the {\em Margulis constant} for $X$. 

For  a nondegenerate arc $I\subset T$  let $\Ga_I$ denote the stabilizer of
$I$ in $\Ga$. Let $\Ga_I'\subset \Ga_I$ be the commutator subgroup. Exhaust $\Ga_I'$
by an increasing sequence of finitely--generated subgroups $\La_n\subset \Ga_I'$.

\begin{lem}\label{elementary}
For each $n$ and $\om$--all $i$, the group $\rho_i(\La_n)$ is elementary.
\end{lem}
\proof The arc $I$ corresponds to a sequence of geodesic arcs
$I_i\subset X$. Let $m_i\in I_i$ be the midpoint.
Let $h_1,...,h_l$ be generators of $\La_n$. Since each $h_j$
is a product of commutators of elements of $\Ga_I$, the arguments of
the proof of Theorem 10.24 in \cite{Kapovich2000} imply that $\rho_i(h_j)$ moves $m_i$
by $\le \mu$ for $\om$--all $i$. Therefore, by Kazhdan--Margulis lemma, the group
$\rho_i(\La_n)$ is elementary. \qed

\medskip
For an elementary subgroup $\La\subset G$,
let $\A(\La)\subset G$ denote the algebraic hull of $\La$ defined in Corollary
\ref{c:hull} and set $A(\La):=\A(\La)\cap \La$.

\medskip
Therefore, each group $\La_{in}:=\rho_i(\La_n)$ contains a canonical
nilpotent subgroup $A_{in}=A(\La_{in})$ of index $\le c$
(where $c$ depends only on $G$). Since $A_{in}$ is canonical, we have
$$
A_{in}\subset A_{i(n+1)}
$$
for each $n$ and $\om$--all $i$. It follows (by taking the $\om$-ultralimit)
that each $\rho_\om(\La_n)$ contains a canonical nilpotent subgroup $A_n$ of index
$\le c$. Thus, the nilpotent subgroup
$$
A:=\bigcup_n A_n
$$
has index $\le c$ in $\rho_\om(\La)$. Therefore, the group $\rho_\om(\Ga_I')$ is
virtually nilpotent. Hence, the group $\rho_\om(\Ga_I)$ fits into the short exact
sequence
$$
1\to \rho_\om(\Ga_I') \to \rho_\om(\Ga_I) \to B \to 1
$$
where $B$ is abelian. Since amenability is stable under group extensions with
amenable kernel and quotient, the group
$\rho_\om(\Ga_I)$ is (elementary) amenable. We proved,
therefore, that $\bar\Ga\acts T$ is {\em  small}. \qed

\begin{rem}
The above argument also works for sequences of group actions on negatively pinched 
Hadamard manifolds of fixed dimensions with fixed pinching constants. 
\end{rem}

The following two propositions describe, to a certain degree, the kernel of the action $\Ga\acts T$.

\begin{prop}\label{ker}
Suppose that each $\Ga_i$ is nonelementary and does not preserve a proper symmetric 
subspace in $X$. Then for every $g\in Ker(\rho_\om)$, for $\om$-all $i$ we have
$$
g\in Ker(\rho_i). 
$$  
\end{prop}
\proof We conjugate the representations $\rho_i$ so that $o_i=o$ for all $i$. 
We will need 

\begin{lem}\label{W}
For every $g\in Ker(\rho_\om)$, we have
$$
\ulim \rho_i(g)=1\in G. 
$$
\end{lem}
\proof For $g\in Ker(\rho_\om)$ set $\ga_i:= \rho_i(g)$. 
Set $R_i:=b(\rho_i)$ and let $B_{R_i}(o)$ be the $R_i$-ball centered at $o$. 

Then we obtain
$$
\ulim \frac{d(x, \ga_i(x))}{R_i}=0, \quad \forall x\in B_{R_i}(o). 
$$
Therefore, there exists $r_i$ so that:
$$
\ulim \frac{r_i}{R_i}\in (0, \infty) 
$$
and for each geodesic segment $\si\subset B_{R_i}(o)$ we have 
$$
dist(\si\cap B_{r_i}(o), \ga_i(\si)\cap B_{r_i}(o))\le \eps_i,   
$$
where 
$$
\ulim \eps_i =0,
$$
and $dist$ stands for the Hausdorff distance. See Lemma 3.10 in \cite{Kapovich2000}. 

By applying this to geodesic segments $\si, \tau\subset B_{R_i}(o)$ which pass through a 
given point $p\in B_1(o)$ and are orthogonal to each other, we conclude that 
$$
\ulim d(\ga_i(p), p)=0. 
$$
Therefore,
$$
\ulim \rho_i(g)=1. \qed 
$$

\medskip 
Let $g_1,...,g_m$ be the generators of $\Ga$. 
Suppose that the assertion of the Proposition fails. Take $\ga\in \ker(\rho_\om)$ 
so that for $\om$-all $i$, $\rho_i(\ga)\ne 1$. By Lemma \ref{W}, 
$$
\ulim ord(\rho_i(\ga))=\infty
$$ 
where $ord$ stands for the order of an element of $G$. 

We now repeat the arguments of the proof of Lemma \ref{kernel}. Let $g\in K=\ker(\rho_\om)$. 
We find finitely-generated subgroups
$$
K_1\subset K_2\subset ... \subset K,
$$
so that $g\in K_1$ and 
$$
g_j K_l g_j^{-1} \subset K_{l+1}, j=1,...,m. 
$$
Lemma \ref{W} implies that for each $l$, 
$\rho_i(K_l)$ is an elementary subgroup of $G$ for $\om$-all $i$. 
Set $\A_{il}:= \A(\rho_i(K_l))$. 
As in the proof of Lemma \ref{kernel}, the group $\A_{il}$ is nontrivial for 
$\om$-all $i$ and each $l$.  

Then, for $\om$-all $i$ there exist $l$ so that for every $j=1,...,m$ we have
$$
\rho_i(g_j) \A_{il} \rho_i(g_j) ^{-1}=\A_{il}. 
$$
Therefore, either $\Ga_i$ is elementary or preserves a proper symmetric subspace in $X$ 
(fixed by $\A_{il}$). In either case, we obtain a contradiction with the assumptions of 
Proposition \ref{ker}. \qed 

\medskip
The tree $T$ contains a unique subtree  $T_{min}$ which is the smallest $\Ga$--invariant 
subtree, see e.g. \cite{Kapovich2000}. The kernel $K$ of the action $\Ga\acts T_{min}$ is, a 
priori, larger than the kernel of $\Ga\acts T$. 

\begin{prop}\label{kern}
Suppose that the tree $T_{min}$ is not a line and the hypothesis of Proposition \ref{ker} hold. 
Then for every $g\in K$, for $\om$-all $i$, we have 
$g\in Ker(\rho_i)$.   
\end{prop}
\proof Since $T_{min}$ is not a line, it contains a nondegenerate triangle 
$x_\om y_\om z_\om\subset T_{min}$. The vertices $x_\om, 
y_\om, z_\om$ of this triangle are represented by sequences $(x_i), (y_i), (z_i)$ in $X$. 
Let $m_i\in \ol{x_i y_i}$ be a point within distance $\le \del$ from the 
other two sides of the triangle $x_i y_i z_i$, where $\del$ is the hyperbolicity constant of $X$.  
For $g\in Ker(\rho_\om)$ set $\ga_i:=\rho_i(g)$. 

\begin{lem}\label{G}
$$
\ulim d(\ga_i(m_i), m_i)=0.
$$
\end{lem}
\proof Our argument is similar to that of the proof of Lemma \ref{W}. 
We again set $R_i:=b(\rho_i)$; then  
$$
\ulim \frac{d(x_i, \ga_i(x_i))}{R_i}=0, \quad \ulim \frac{d(y_i, \ga_i(y_i))}{R_i}=0,  
$$
$$
\ulim \frac{d(z_i, \ga_i(z_i))}{R_i}=0.  
$$
As in the proof of Lemma \ref{W}, the segment $\ol{x_i y_i}$ will 
contain a subsegment $\si_i:=\ol{x_i' y_i'}$ of length $r_i$ so that 
$m_i\in \ol{x_i' y_i'}$,  
$$
\ulim \frac{r_i}{R_i}\in (0, \infty),  
$$
and 
$$
\ulim d(\ga_i(x_i'), \si_i)=\ulim d(\ga_i(y_i'), \si_i)=0.
$$

Define points  $p_i, q_i\in \si_i$ nearest to $z_i, \ga_i(z_i)$ respectively. 
Then
$$
\ulim d(q_i, \ga_i(p_i))=0. 
$$

\begin{figure}[tbh]
\begin{center}
\input{f.pstex_t}
\end{center}
\caption{} \label{f.fig}
\end{figure}

Suppose that the isometries $\ga_i$ {\em shear} along the segments $\si_i$,  i.e. 
$$
\ulim d(\ga_i(m_i), m_i)\ne 0. 
$$
Then   
$$
\ulim d(\ga_i(p_i), p_i)= \ulim d(p_i, q_i)=\ulim d(\ga_i(m_i), m_i)\ne 0. 
$$
 Since 
$$
\ulim \frac{d(z_i, p_i)}{R_i}=d(z_\om, p_\om)\ne 0, \quad \ulim \frac{d(z_i, \ga_i(z_i))}{R_i}= 0,
$$
it follows that there exists a point $w_i\in \ol{z_i \ga_i(z_i)}$   
within distance $\le 2\del$ from both 
$$
\ol{z_i p_i}, \quad \ol{q_i \ga_i(z_i)}. 
$$ 
See Figure \ref{f.fig}. Since 
$$
\ulim \frac{d(z_i, \ga_i(z_i))}{R_i}=0, \quad \ulim \frac{d(z_i, p_i)}{R_i}\ne 0, 
$$
we obtain 
$$
\ulim d(w_i, \si_i)=\infty. 
$$

Take the shortest segments $\rho_i, \tau_i$ from $w_i$ to  $\ol{z_i p_i}, \ol{q_i \ga_i(z_i)}$. 
The nearest-point projection to $\si_i$ sends $\rho_i\cup \tau_i$ onto $\ol{p_i q_i}$. 
However, this projection is exponentially contracting and 
$\ulim d(w_i, \si_i)=\infty$. 
This contradicts the assumption that
$$
\ulim d(p_i, q_i)\ne 0. 
$$
Therefore, 
$$
\ulim d(\ga_i(m_i), m_i)= 0. \qed
$$
 
Given $g\in K=\ker(\Ga\to \Isom(T_{min}))$, 
we define the finitely-generated subgroups $K_l\subset K$ 
in the same fashion it was done in the proof of Proposition \ref{ker}. By Lemma \ref{G}, it 
follows that for every generator $h\in K_l$ and $\eta_i:=\rho_i(h)$, we have 
$$
\ulim d(\eta_i(m_i), m_i)=0. 
$$ 
Therefore, by Kazhdan--Margulis lemma, for each $l$ and $\om$-all $i$, the group 
$\rho_i(K_l)$ is elementary. Now, the arguments from the proof of Proposition \ref{ker} 
go through and we obtain $\rho_i(g)=1$ for $\om$-all $i$. \qed

\section{Semistability}\label{semistability}

The purpose of this section is to weaken the notion of {\em stability}
used in the Rips' theory, so that the Rips Machine still applies. We recall

\begin{dfn}
Let $\Ga\acts T$ be an isometric group action on a tree.
A nondegenerate arc $I\subset T$ is called {\em stable} if for every decreasing
sequence of nondegenerate subarcs
$$
I\supset I_1\supset I_2 \supset ...$$
the corresponding sequence of stabilizers
$$
\Ga_I\subset \Ga_{I_1}\subset \Ga_{I_2} \subset ...$$
is eventually constant. The action  $\Ga\acts T$ is called {\em stable} if
every nondegenerate arc $J\subset T$ contains a stable subarc.
\end{dfn}

M.~Dunwoody \cite{Dunwoody(1999)} constructed example of a small
but unstable action of a finitely generated group $\Ga$ on a tree. To remedy this, we
introduce the following modification of stability, 
adapted to the case of actions whose image on $\Isom(T)$ is small:

\begin{dfn}\label{sst}
Suppose that we are given an isometric action of a group on a tree $\rho: \Ga\acts T$.
We say that this action  is {\em semistable} if it satisfies
the following property:

For every nondegenerate arc $I\subset T$ and its stabilizer
$\Ga_I\subset \Ga$, there exists a canonical amenable subgroup
$\A(\Ga_I)\subset \Isom(T)$ so that:

1. If $I\supset J$ then $\A(\Ga_I)\subset \A(\Ga_J)$.

2. $A(\Ga):= \A(\Ga_I)\cap \rho(\Ga_I)$ has index $\le c<\infty$ in $\rho(\Ga_I)$,
where $c=c_T$ is a uniform constant.

3. If $\al\in \Ga$ is such that $\al \A(\Ga_I) \al^{-1}\subset \A(\Ga_I)$, then
$\al \A(\Ga_I) \al^{-1}= \A(\Ga_I)$.

4. For every nondegenerate arc $J\subset T$, there exists a nondegenerate subarc
$I\subset J$ so that the following holds:

If $I\supset I_1\supset I_2 \supset ... $ is a decreasing sequence of nondegenerate
arcs, then the sequence of groups
$$
\A(\Ga_1)\subset  \A(\Ga_2)\subset ...
$$
is eventually constant.
\end{dfn}

We say that an arc $I\subset T$ is {\em stabilized} (with respect to the
action of $\Ga$) if for every nondegenerate subarc $J\subset I$, we have
$$
\A(\Ga_I)=\A(\Ga_J).
$$
We let $\A_I$ denote $\A(\Ga_I)$ in this case.

Note that every semistable action is automatically small, since a finite index
extension of an amenable group is also amenable.

\medskip
It is easy to classify the possible amenable groups $\A\subset \Isom(T)$:

1. $\A$ is parabolic, i.e., it fixes a point in $\geo T$ and does
not fix any other points in $T\cup =\geo T$.

2. $\A$ is hyperbolic, i.e., it has a unique invariant geodesic
$T_\A\subset T$ and contains a
nontrivial translation along this geodesic.

3. $\A$ is elliptic, i.e., it fixes a nonempty subtree $T_\A\subset T$.

\medskip

We now give examples of semistable actions.

\smallskip
{\bf Example 1.} Consider $\rho: \Ga\to \Isom(T)$, so that  
the action of the image group $\bar\Ga=\rho(\Ga)$ on $T$ is small and stable. 
Then $\Ga\acts T$ is also semistable: take $\A(\Ga_I):=\rho(\Ga_I)$.

\medskip
{\bf Example 2.} Let $\F$ be
a nonarchimedean valued field of zero characteristic and cardinality continuum
and $G=\ul{G}(\F)$ be a group of rank 1. We then consider the Bruhat--Tits tree $T$
associated with the group $G$. The quotient group $G/Z(G)$ acts faithfully on $T$,
where $Z(G)$ is the center of $G$. Let $\Ga\subset G/Z(G)\subset \Isom(T)$ be a
 subgroup so that the associated action $\Ga\acts T$ is small.

Given an amenable subgroup $\La\subset \Ga$, consider its lift $\t{\La}\subset G$,
which is still an amenable subgroup. Let $\A(\La)\subset G/Z(G)$ denote the
projection of the hull $\A(\t\La)\subset G$,  defined
in Corollary \ref{c:hull}. It is immediate that $\A(\La)$
satisfies Properties 1--3 of Definition \ref{sst}. Consider
Property 4.

For the amenable groups $\La=\Ga_I$,
the algebraic hulls $\A(\Ga_I)$ are Zariski connected algebraic subgroups of $G$.
Since the dimensions of the groups in the sequence
$$
\A(\Ga_I)\subset \A(\Ga_{I_1})\subset \A(\Ga_{I_2}) ...
$$
are eventually constant,  this sequence is eventually constant as well.

\medskip
{\bf Example 3.} Let $\rho_\om: \Ga\to Isom(T)$ be a group action on a tree associated
with a divergent sequence of representations $\rho_i: \Ga\to \Isom(X)$, where $X$
is a negatively curved symmetric space. The  asymptotic cone  $T=T_\om$ of $X$
associated with this sequence is a metric tree. According to
\cite{Chiswell, KSTT, Thornton}, the asymptotic cone $T$ is the Bruhat-Tits tree
associated with a group $\ul{G}(\F)$, where $\F$ is a certain nonarchimedean
valued complete field of cardinality continuum and characteristic zero.
Moreover, the group $\Ga$ maps to $\Isom(T)$ via a homomorphism
$$
\rho_\om: \Ga \to \ul{G}(\F)\subset \Isom(T).
$$

\begin{rem}
The field $\F$ is a subfield
of the field of nonstandard reals, which is
the ultraproduct
$$
\R_*=\prod_{i\in \N} \R  /\om.
$$
The choice of the subfield and valuation depends on $\om$ and on the divergent
sequence $b(\rho_i)$.

In the case $X=\H^3$ and $\ul{G}=SL(2)$, we can use algebraically closed field $\F$,
which is a subfield of the ultraproduct
$$
\C_*=\prod_{i\in \N} \C / \om.
$$
\end{rem}

Therefore, for each amenable subgroup $\La=\rho_\om(\Ga_I)$,
we can define the algebraic hull $\A(\La)$ using Corollary \ref{c:hull}
(see Example 2 above).
In case $X=\H^3$ and $\F$ algebraically closed, we can use Theorem \ref{hull}, or,
rather, the example which appears in the beginning of the proof. In particular, by Example 2, 
it follows that the action $\rho_\om: \Ga\acts T$ is semistable.

\begin{cor}
\label{sstability}
Part 2 of Theorem \ref{main2} holds.
\end{cor}

\medskip
{\bf Implications of semistability.}
We now assume that we are given a semistable action $\Ga\acts T$ and the corresponding 
action $\bar\Ga\acts T$ of the image of $\Ga$ in $\Isom(T)$.

Let $I\subset T$ be a stabilized arc and $\al\in \Ga$ be an axial isometry of $T$,
whose axis contains $I$, and so that
$$
J=I\cap \al(I)
$$
is nondegenerate. Then
$$
\al \Ga_I \al^{-1} \subset \Ga_J.
$$
Since $I$ is stabilized, follows that $\al \A_I \al^{-1} \subset \A_J=\A_I$. Thus  
$\al \A_I \al^{-1} =\A_I$ (see Part 3 of Definition \ref{sst}).
Suppose that we are given two elements $\al, \be \in \Ga$ as above, so that
\begin{equation}
\label{id}
\al \A_I \al^{-1} =\A_I, \quad \be \A_I \be^{-1} =\A_I.
\end{equation}

{\bf Case 1.} $\A_I$ is parabolic. Then the equalities (\ref{id}) imply that
$\al, \be$ both fix the unique fixed point at infinity of
the group $\A_I$. Since the action $\bar\Ga\acts T$ is small, 
it follows that the group $\rho\<\al, \be\>$ generated by 
$\rho(\al), \rho(\be)$ is virtually solvable, see \cite[\S 10.5]{Kapovich2000}.

{\bf Case 2.}  $\A_I$ is hyperbolic. Then the equalities (\ref{id})
imply that $\al, \be$ preserve the unique invariant geodesic of the group $\A_I$.
Hence, the commutator subgroup of $\<\al, \be\>$ fixes this geodesic pointwise.
It again follows that $\rho\<\al, \be\>$ is virtually solvable.

{\bf Case 3.} $\A_I$ is elliptic. Let $T'\subset T$ denote the subtree fixed by $\A_I$.
Then $T'$ is invariant under both $\al$ and $\be$. The restrictions of these isometries
to $T'$ remain axial.

Recall that
$$
|\rho(\Ga_J): A(\Ga_J)|\le c_T
$$
for every arc $J$.

\begin{ass}\label{in}
We now assume in addition that $n$ is a natural number so that
$$
\frac{\length(A_\al\cap A_\be)}{\ell(\al)+\ell(\be)}\ge 2n> 2c_T.
$$
\end{ass}

Under this assumption, for each $i=1,..,n$, $[\al^i, \be]\in \Ga_J\subset \bar\Ga_J$
for some nondegenerate
subinterval $J\subset I$. Moreover, there exist $m\ne n$ so that we have the equality of the
cosets
$$
[\al^m, \be] \A_I= [\al^n, \be] \A_I.
$$
Since $\A_I$ fixes $T'$ pointwise, it follows that
$$
[\al^m, \be]|_{T'}= [\al^n, \be]|_{T'}.
$$
Hence,
$$
[\al^{m-n}|_{T'}, \be|_{T'}] =1.
$$

Since $\al^{m-n}|_{T'}, \be|_{T'}$ are commuting nontrivial axial elements, they
have to have common axis. Therefore, $\al, \be$ also have common axis.
Now, analogously to the Case 2, it follows that $\rho\<\al, \be\>$ is virtually solvable.

\medskip
We conclude that in each case (provided that the Assumption \ref{in} holds in the
elliptic case), we have

\begin{prop}\label{amenable}
The group $\rho\<\al, \be\>$ is amenable.
\end{prop}

\section{Generalization of the Rips theory\\ to the semistable case}

In this section we will finish the proof of Theorem \ref{main2} by verifying Part 3. 

Suppose that we are given a semistable nontrivial action
$$
\rho: \Ga\to \Isom(T),$$
of a finitely-presented $\Ga$ on a tree $T$. 
Then one can apply the arguments of the Rips Theory (see \cite{Bestvina-Feighn(1995)}
or \cite[Chapter 12]{Kapovich2000}) to the action $\Ga\acts T$. Note that the only
place the {\em stability} condition is used in the proof of the Rips theorem, is
the analysis of the {\em axial} pure band complex $C$, see e.g.
\cite[Proposition 12.69]{Kapovich2000}.

In this case one deals with pairs of axial isometries
$\al, \be \in \Ga$, so that the ratio
$$
\frac{\length(A_\al\cap A_\be)}{\ell(\al)+\ell(\be)}
$$
can be taken as large as one wishes. Therefore, one can choose this ratio to satisfy
the Assumption \ref{in} as above. The conclusion of the Rips Theory in the Axial
case is then that the action of the fundamental group $\pi_1(C)$ of the component
$C$ (which is a subgroup of $\Ga$) on the tree $T$ has an invariant geodesic.
It then deduced that the action of  $\pi_1(C)$ factors through action of a
solvable group.

In our case, Proposition \ref{amenable} implies that the action  $\pi_1(C)\acts T$
either has an invariant geodesic or is parabolic; in either case,
it factors through action of an amenable group.

Therefore,  repeating verbatim the proof of
Theorem 12.72 in \cite{Kapovich2000}, we obtain

\begin{thm}\label{rips}
One of the following holds:

1. If the action $\Ga\acts T$ is not {\em pure} then  the group $\Ga$
splits nontrivially as $\Ga=\Ga_1*_E \Ga_2$ or $\Ga=\Ga_1*_E$, over a subgroup $E$,
which fits into a short exact sequence
$$
1\to K_E\to E \to Q\to 1,
$$
where $K_E$ fixes a nondegenerate arc in $T$ and $Q$ is either finite or cyclic.
Moreover, the group $E$ fixes a point in $T$ and 
the groups $\Ga_1, \Ga_2, E$ are finitely generated. 

2. If the action is pure then $G$ belongs to one of the following types:

(a) Surface type.

(b) Axial type.

(c) Thin type.

In either case, $\Ga$ splits nontrivially
as $\Ga=\Ga_1*_E \Ga_2$ or $\Ga=\Ga_1*_E$, over a subgroup $E$, which fits into a
short exact sequence
$$
1\to K_E\to E \to Q\to 1,
$$
where $K_E$ fixes a nondegenerate arc in $T$ and $Q$ is abelian. 
The groups $\Ga_1, \Ga_2, E$ are finitely generated.
\end{thm}

Therefore, the image (in $\Isom(T)$) of the edge subgroup of $\Ga$ is amenable.

We now assume that the action $\Ga\acts T$ arises from a divergent
sequence of discrete but not necessarily
faithful representations
$$
\rho_i: \Ga \to \Isom(X)
$$
where $X$ is a negatively curved symmetric space. Then we obtain $\Ga\acts T$,
where $T$ is an asymptotic cone of $X$, which can be realized
as the Bruhat-Tits tree of a rank 1 algebraic group $\ul{G}(F)$. Thus we obtain
a homomorphism $\rho_\om: \Ga\to \bar\Ga\subset \ul{G}(F)\subset \Isom(T)$.
Then, according to Section \ref{small}, the action $\bar\Ga\acts T$ is small.
According to Section  \ref{semistability}, this action is also
semistable. Therefore, Theorem \ref{rips} applies and we obtain:

\begin{cor}
The group $\Ga$ splits  nontrivially
as $\Ga=\Ga_1*_E \Ga_2$ or $\Ga=\Ga_1*_E$, over a subgroup $E$, so
that $\rho(E)$ is amenable. The groups $\Ga_1, \Ga_2, E$ are finitely generated. 
\end{cor}

\begin{rem}
M.~Dunwoody \cite{Dunwoody(2006)} proved another version of Rips Theorem
in the case of {\em slender} faithful actions of finitely-presented groups
on trees without the stability hypothesis.
However his main theorem only yields a splitting of $\Ga$ where
each edge group is either slender or fixes a point in the tree. This
is not enough to guarantee amenability of the edge groups in the resulting
decomposition. Moreover, it appears that the arc stabilizers $\Ga_I$
for group actions on trees associated with divergent sequences of discrete
representations, need not be slender. For instance, it seems that they can
contain infinitely generated abelian subgroups.
\end{rem}

Since $\rho(E)\subset \ul{G}(F)$, it follows that this subgroup is virtually solvable. 
By combining the above results, we obtain Theorem \ref{main2}.

\newpage

\bibliography{../lit}

\begin{thebibliography}{10}

\bibitem{Ballmann-Gromov-Schroeder(1985)}
{\sc W.~Ballmann, M.~Gromov, and V.~Schroeder}, {\em Manifolds of Nonpositive
  Curvature}, Progress in Math., vol. 61, Birkh{\"a}user, 1985.

\bibitem{Belegradek1998}
{\sc I.~Belegradek}, {\em Intersections in hyperbolic manifolds}, Geometry and
  Topology (electronic), 2 (1998), pp.~117--144.

\bibitem{Belegradek-Kapovitch}
{\sc I.~Belegradek and V.~Kapovitch}, {\em Classification of negatively pinched
  manifolds with amenable fundamental groups}, Acta Math., 196 (2006),
  pp.~229--260.

\bibitem{Bestvina(1988)}
{\sc M.~Bestvina}, {\em Degenerations of hyperbolic space}, Duke Math. Journal,
  56 (1988), pp.~143-- 161.

\bibitem{Bestvina-Feighn(1995)}
{\sc M.~Bestvina and M.~Feighn}, {\em Stable actions of groups on real trees},
  Inventiones Mathematicae, 121 (1995), pp.~287--321.

\bibitem{Chiswell}
{\sc I.~M. Chiswell}, {\em Nonstandard analysis and the {M}organ-{S}halen
  compactification}, Quart. J. Math. Oxford Ser. (2), 42 (1991), pp.~257--270.

\bibitem{Chuckrow}
{\sc V.~Chuckrow}, {\em Schottky groups with applications to {K}leinian
  groups}, Ann. of Math., 88 (1968), pp.~47--61.

\bibitem{Dunwoody(1999)}
{\sc M.~J. Dunwoody}, {\em A small unstable action on a tree}, Math. Res.
  Lett., 6 (1999), pp.~697--710.

\bibitem{Dunwoody(2006)}
\leavevmode\vrule height 2pt depth -1.6pt width 23pt, {\em Groups acting on
  real trees}.
\newblock Preprint, 2006.

\bibitem{GLP2}
{\sc D.~Gaboriau, G.~Levitt, and F.~Paulin}, {\em Pseudogroups of isometries of
  {${\mathbb R}$} and constructions of {${\mathbb R}$}-trees}, Ergodic Theory
  and Dynamical Systems, 15 (1995), pp.~633--652.

\bibitem{Guirardel(2006)}
{\sc V.~Guirardel}, {\em Actions of finitely generated groups on $\r$-trees}.
\newblock Preprint, math/0607295, 2006.

\bibitem{Jorgensen(1976)}
{\sc T.~Jorgensen}, {\em On discrete groups of {M}\"obius transformations},
  Amer. J. Math., 98 (1976), pp.~739--749.

\bibitem{Jorgensen(1982)}
{\sc T.~Jorgensen and P.~Klein}, {\em Algebraic convergence of finitely
  generated {K}leinian groups}, Quart. J. Math. Oxford, 33 (1982),
  pp.~325--332.

\bibitem{Kapovich2000}
{\sc M.~Kapovich}, {\em Hyperbolic manifolds and discrete groups}, Birkh\"auser
  Boston Inc., Boston, MA, 2001.

\bibitem{Kapovich-Leeb(1995)}
{\sc M.~Kapovich and B.~Leeb}, {\em On asymptotic cones and quasi-isometry
  classes of fundamental groups of $3$-manifolds}, Journal of Geometric and
  Functional Analysis, 5 (1995), pp.~582--603.

\bibitem{KSTT}
{\sc L.~Kramer, S.~Shelah, K.~Tent, and S.~Thomas}, {\em Asymptotic cones of
  finitely presented groups}, Adv. Math., 193 (2005), pp.~142--173.

\bibitem{Martin(1989a)}
{\sc G.~Martin}, {\em On discrete {M}obius groups in all dimensions: A
  generalization of {J}orgensen's inequality}, Acta Math., 163 (1989),
  pp.~253--289.

\bibitem{Martin(1993)}
\leavevmode\vrule height 2pt depth -1.6pt width 23pt, {\em On discrete isometry
  groups of negative curvature}, Pacific J. Math., 160 (1993), pp.~109-- 128.

\bibitem{Morgan(1986)}
{\sc J.~Morgan}, {\em Group actions on trees and the compactification of the
  space of classes of ${S}{O}(n, 1)$ representations}, Topology, 25 (1986),
  pp.~1--33.

\bibitem{Morgan-Shalen(1984)}
{\sc J.~Morgan and P.~Shalen}, {\em Valuations, trees and degenerations of
  hyperbolic structures {I}}, Ann. of Math., 120 (1984), pp.~401--476.

\bibitem{Paulin(1988)}
{\sc F.~Paulin}, {\em Topologie de {G}romov equivariant, structures
  hyperboliques et arbres reels}, Inventiones Mathematicae, 94 (1988), pp.~53--
  80.

\bibitem{Paulin(1996)}
\leavevmode\vrule height 2pt depth -1.6pt width 23pt, {\em Actions de groupes
  sur les arbres}, S\'eminaire Bourbaki,  (1997), pp.~97--137.

\bibitem{Poincare}
{\sc H.~Poincar\'e}, {\em On the groups of linear equations}, in ``Papers on
  {F}uchsian {F}unctions'', Springer Verlag, 1985, pp.~357--483.

\bibitem{Raghunathan}
{\sc M.~Raghunathan}, {\em Discrete subgroups of Lie groups}, Springer, 1972.

\bibitem{Rips-Sela(1994)}
{\sc E.~Rips and Z.~Sela}, {\em Structure and rigidity in hyperbolic groups,
  {I}}, Journal of Geometric and Functional Analysis, 4 (1994), pp.~337--371.

\bibitem{Thornton}
{\sc B.~Thornton}, {\em Asymptotic cones of symmetric spaces}.
\newblock Ph.{D}. {T}hesis, {U}niversity of {U}tah, 2002.

\bibitem{Wang-Yang}
{\sc X.~Wang and W.~Yang}, {\em Discreteness criteria of {M}\"obius groups of
  high dimensions and convergence theorems of {K}leinian groups}, Adv. Math.,
  159 (2001), pp.~68--82.

\bibitem{Wielenberg}
{\sc N.~Wielenberg}, {\em Discrete {M}oebius groups: fundamental polyhedra and
  convergence}, Amer. Journ. Math., 99 (1977), pp.~861--878.

\end{thebibliography}
\bibliographystyle{siam}

\noindent Department of Mathematics, 1 Shields Ave.,\\
University of California, Davis, CA 95616, USA\\
kapovich@math.ucdavis.edu

\end{document}